\documentclass[graybox]{svmult}

\usepackage{type1cm}        % activate if the above 3 fonts are
\usepackage{makeidx}         % allows index generation
\usepackage{graphicx}        % standard LaTeX graphics tool
\usepackage{multicol}        % used for the two-column index
\usepackage[bottom]{footmisc}% places footnotes at page bottom

\usepackage{newtxtext}       % 
\usepackage{newtxmath}       % selects Times Roman as basic font

\makeindex             % used for the subject index
\usepackage{tikz} 
\usepackage{pgfplots}
\newcommand{\matr}[1]{\mathbf{#1}}
\newcommand{\Qmat}{\matr{Q}}
\newcommand{\QDmat}{\matr{Q}_\Delta}

\newcommand{\vect}[1]{\boldsymbol{#1}}
\newcommand{\dt}{\Delta t}

\begin{document}

\title*{The Parallel Full Approximation Scheme in Space and Time for a Parabolic Finite Element Problem}
\titlerunning{Coupling PFASST with Finite Elements in Space}
%\title*{Contribution Title}
\author{Oliver Sander \and Ruth Sch\"obel \and Robert Speck}
\institute{Oliver Sander \at Faculty of Mathematics\\ TU Dresden, Germany \email{oliver.sander@tu-dresden.de}
    \and Ruth Sch\"obel \at J\"ulich Supercomputing Centre\\ Forschungszentrum J\"ulich GmbH, Germany, \email{r.schoebel@fz-juelich.de}
	\and Robert Speck \at J\"ulich Supercomputing Centre\\ Forschungszentrum J\"ulich GmbH, Germany, \email{r.speck@fz-juelich.de}	
}
\maketitle

\abstract*{}

\abstract{
The parallel full approximation scheme in space and time (PFASST) is a parallel-in-time integrator that allows to integrate multiple time-steps simultaneously.
It has been shown to extend scaling limits of spatial parallelization strategies when coupled with finite differences, spectral discretizations, or particle methods.
In this paper we show how to use PFASST together with a finite element discretization in space.
While seemingly straightforward, the appearance of the mass matrix and the need to restrict iterates as well as residuals in space makes this task slightly more intricate. 
We derive the PFASST algorithm with mass matrices and appropriate prolongation and restriction operators and show numerically that PFASST can, after some initial iterations, gain two orders of accuracy per iteration.
}

\section{Introduction}

The parallel full approximation scheme in space and time  (PFASST, \cite{EmmettMinion2012}) can integrate multiple time-steps simultaneously by using inner iterations of spectral deferred corrections (SDC, \cite{DuttEtAl2000}) on a space-time hierarchy. 
It mimics a full approximation scheme (FAS, \cite{trottenberg}) for a sequence of coupled collocation problems. 
For the simulation of space--time dependent problems, PFASST has been used in combination with finite differences, e.g., in~\cite{minion2015interweaving,schobel2019pfasst}, but also in connection with particle simulations \cite{speck_massively_2012} and spectral methods~\cite{gotschel2017parallel}.
In this work we combine PFASST with a finite element discretization in space.
Using a simple, nonlinear reaction--diffusion equation, we will derive the discretized, ``composite collocation problem'' PFASST aims to solve in parallel and show the correct handling of the mass matrix. 
There exist two different ways to write down the composite collocation problem with a non-trivial mass matrix and we will demonstrate that we can avoid inversion of the mass matrix with the added benefit of a better order of accuracy in time per PFASST iteration.
The choice of restriction and prolongation in space plays a major role and we will show the correct formulation and placement of those.
Both mass matrix handling and choice of transfer operators mark key differences to using standard finite differences in space, both in terms of theoretical formulation and simulation results.
Using a concrete example, we numerically test the order of accuracy per iteration of PFASST and compare it with SDC.

\section{PFASST and finite elements in space}\label{sec_spatial}
We consider the reaction--diffusion equation
\begin{align}
   v_t(x,t) &= \Delta v(x,t) + \text{g}(v(x,t)),  &&  x \in \Omega, t \in [t_0,T], \label{s1} \\
   v(x,t) &= 0, && x \in \partial \Omega, \nonumber 
  \end{align}  
with suitable initial conditions for $t = t_0$ and $g: \mathbb{R} \rightarrow \mathbb{R}$ continuously differentiable. 
Here $\Omega \subset \mathbb{R}$ is a polyhedral domain with boundary $\partial \Omega$, and $\Delta$ denotes the Laplace operator.

\subsection{Finite element discretization in space}\label{mol_space}
We define test functions $\phi^h$ in a finite-dimensional space $\matr V^h \subset \matr H^1_0(\Omega)$, 
multiply \eqref{s1} by these test functions, and integrate by parts. Thus, $v^h(\cdot, t) \in \matr V^h$ is given by 
\begin{align}
\int_{\Omega} \phi^h v^h_t \ dx =  - \int_{\Omega} \nabla \phi^h \nabla v^h \ dx + \int_\Omega \phi^h \text{g}(v^h) \ dx \quad \forall \phi^h \in \matr V^h  \label{eq:disc}. 
\end{align}
We choose a basis $\varphi_1, \dots, \varphi_N$ of $\matr V^h$ 
and approximate $\text{g}(v^h)$ by an element of $\matr V^h$ and express $v^h$ and $\text{g}(v^h)$ as
\begin{align}
v^{h}(x,t) =   \sum\limits_{i=1}^N v_i(t) \varphi_i(x), \qquad \text{g}(v^{h})(x,t) \approx   \sum\limits_{i=1}^N \text{g}(v_i(t)) \varphi_i(x), \label{eq:ling}
\end{align}
where the coefficients $v_i(t), \ i=1, \dots, N$, are time-dependent functions.
Inserting~\eqref{eq:ling} into equation \eqref{eq:disc} yields 
\begin{align}
 \matr{M} u_t = - \matr{A} u + \matr M g(u) \eqqcolon f(u). \label{final_rothe}
\end{align}
Here, $u \coloneqq (v_1, \dots, v_N)$ is a vector holding the coefficients $v_i$, and $g: \mathbb{R}^N \rightarrow \mathbb{R}^N$, $g \coloneqq (\text{g}(v_{1}), \dots,  \text{g}(v_{N}))^T$.
The matrix $\matr M \in \mathbb{R}^{N \times N}$ is the mass matrix and $\matr A \in \mathbb{R}^{N \times N}$  the stiffness matrix 
 \begin{align*}
  \matr{M}_{ij} \coloneqq \int_{\Omega_i} \varphi_i \varphi_j \ dx, \qquad  	 \matr{A}_{ij} \coloneqq  \int_{\Omega_i} \nabla \varphi_i \nabla \varphi_j \ dx. 
\end{align*}
%respectively.

\subsection{The collocation problem and SDC}

For the temporal discretization, we decompose the interval $[t_0,T]$ into time-steps $t_0 < t_1 < \dots < t_L =T$, $L\in\mathbb{N}$.
For one time-step $[t_l,t_{l+1}]$, the Picard formulation of~\eqref{final_rothe} is
\begin{align}\label{eq:ODE2}
 \matr M  u(t) = \matr M u_{l,0} + \int_{t_l}^t f(u(s)) \ ds,\qquad t\in[t_l,t_{l+1}],
\end{align}
where $u_{l,0}\coloneqq u(t_l)$.
To approximate the integral we use a spectral quadrature rule on $[t_l, t_{l+1}]$
with $M$ quadrature nodes $\tau_{l,1},...,\tau_{l,M}$ such that $t_l < \tau_{l,1} < ... < \tau_{l,M} = t_{l+1}$.%(note the Gau\ss-Lobatto property $\tau_{l,M}=t_{l+1}$).

For each of the $M$ nodes we introduce a set of $M$ quadrature weights $q_{m,j} \coloneqq \int_{t_l}^{\tau_{l,m}} L_j(s)\ ds$, $m,j = 1, \dots , M$, where $L_1, \dots, L_M$ are the Lagrange polynomials for the nodes $\tau_{l,1}, ..., \tau_{l,M}$.
We can then approximate the integral in \eqref{eq:ODE2} from $t_l$ to %the nodes 
$\tau_{l,m}$ by
\begin{align*} %\label{integrals}
  \Delta t \sum_{j=1}^M q_{m,j}f(u_{l,j}) \approx \int_{t_l}^{\tau_{l,m}} f(u(s)) \ ds, \quad m= 1, \dots, M,
\end{align*} 
where $\dt \coloneqq t_{l+1}-t_l$ denotes the time-step size. 
Using this in Equation \eqref{eq:ODE2}  the unknown values $u(\tau_{l,1}), \dots, u(\tau_{l,M})$ can be approximated by a solution $u_{l,1}, \dots, u_{l,M} \in \mathbb{R}^N$ of the nonlinear system of equations
\begin{align*}%\label{equations}
  \matr M u_{l,m} = \matr M u_{l,0} + \Delta t \sum_{j=1}^M q_{m,j}f(u_{l,j}) \quad \text{for } m= 1, \dots, M.
\end{align*} 
This is the so called ``collocation problem'', which we can rewrite as  
\begin{align}\label{eq:coll_prob}
     \matr{C}^{\operatorname{coll}}_{\vect f}(\vect{u}_l) \coloneqq \left(\matr{I}_{M} \otimes \matr M  - \dt(\Qmat \otimes \matr I_N)\vect{f} \right)(\vect{u}_l) = (\matr I_M \otimes \matr M ) \vect{u}_{l,0},
\end{align}
where $\matr I_X \in \mathbb{R}^{X \times X}$, $X\in\mathbb{N}$ is the identity matrix, $\otimes$ denotes the Kronecker product, $\vect{u}_l \coloneqq (u_{l,1}, ..., u_{l,M})^T \in \mathbb{R}^{MN}$, $\vect{u}_{l,0} \coloneqq (u_{l,0}, ..., u_{l,0})^T\in\mathbb{R}^{MN}$, $\Qmat \coloneqq (q_{ij})\in\mathbb{R}^{M\times M}$, and the vector function $\vect{f}:\mathbb{R}^{MN} \to \mathbb{R}^{MN}$ is given by 
$
    \vect{f}(\vect{u}_l) \coloneqq  (f(u_{l,1}), ..., f(u_{l,M}))^T.
$

With this matrix notation, spectral deferred corrections can simply be seen as a preconditioned Picard iteration~\cite{HuangEtAl2006,RuprechtSpeck2016}. 
More precisely, for a lower triangular matrix $\QDmat \in \mathbb{R}^{M\times M}$ we define the preconditioner
\begin{align*}%\label{star}
  \matr{P}^{\operatorname{sdc}}_{\vect f}(\vect{u}_l) \coloneqq \left(  \matr{I}_{M} \otimes \matr M - \dt(\QDmat \otimes \matr I_N)\vect{f} \right)(\vect{u}_l).
\end{align*}
Then the preconditioned iteration reads
\begin{align} \label{eq:SDC}
 \matr{P}^{\operatorname{sdc}}_{\vect f}(\vect{u}_l^{k+1}) 
 &= (\matr{P}^{\operatorname{sdc}}_{\vect f} -  \matr{C}^{\operatorname{coll}}_{\vect f})(\vect{u}_l^k) + (\matr I_M \otimes \matr M) \vect{u}_{l,0},\quad k=1, ..., K.
\end{align} 

The properties of $\matr{P}^{\operatorname{sdc}}_{\vect f}$ depend first and foremost on the choice of the matrix $\QDmat$.
For this work, we use the backward Euler approach.
We refer to~\cite{HuangEtAl2006,Weiser2014,RuprechtSpeck2016} for more details on the notation and its relationship to the original description of SDC as in~\cite{DuttEtAl2000}.
The key difference to the usual formulation of the collocation problem~\eqref{eq:coll_prob} and SDC~\eqref{eq:SDC} iteration is the appearance of the mass matrix $\matr{M}$, which for finite difference discretization is just the identity matrix.

\subsection{The composite collocation problem and PFASST}

For $L$ time-steps, the composite collocation problem is
 \begin{align} \label{ccOp}
  \begin{pmatrix}
    \matr{C}^{\operatorname{coll}}_{\vect f} \\
     -\matr{H} & \matr{C}^{\operatorname{coll}}_{\vect f} \\
      & \ddots & \ddots \\
     & & -\matr{H} & \matr{C}^{\operatorname{coll}}_{\vect f}
  \end{pmatrix}
  \begin{pmatrix}
     \vect{u}_{1}\\
     \vect{u}_{2}\\
     \vdots\\
     \vect{u}_{L}
  \end{pmatrix} = 
  \begin{pmatrix}
    (\matr I_M \otimes \matr M) \vect{u}_{0,0} \\
     \vect 0\\
     \vdots\\
     \vect 0
  \end{pmatrix},
 \end{align}
where the matrix $\matr{H}\coloneqq \matr N \otimes \matr M$ with $\matr{N}\coloneqq (\delta_{i,M})_{i,j} \in\mathbb{R}^{M\times M}$ provides the value at the last quadrature node $\tau_{l,M}$ of a time-step $[t_l, t_{l+1}]$ as initial value for the following time-step. 
Defining the global state vector $\vect{u} \coloneqq (\vect{u}_{1}, ..., \vect{u}_{L})^T\in\mathbb{R}^{LMN}$, the vector $\vect{b} \coloneqq ((\matr I_M \otimes \matr M) \vect{u}_{0,0}, \vect 0, ..., \vect 0)^T\in\mathbb{R}^{LMN}$,
and ${\vect F: \mathbb{R}^{LMN} \rightarrow \mathbb{R}^{LMN}}$ with
${\vect F} ({\vect u}) \coloneqq \left( {\vect f} ({\vect u}_{1}), \dots , {\vect f} ({\vect u}_{L})  \right)^T$,
we can write this in the more compact form as $\matr C_{\vect F}(\vect u) = \vect b$, where $\matr C_{\vect F}$ is the lower block-bidiagonal, nonlinear operator on the left of \eqref{ccOp}.
Using the definition~\eqref{eq:coll_prob} of $\matr{C}^{\operatorname{coll}}_{\vect f}$ we write~\eqref{ccOp} as
\begin{align}\label{eq:comp_coll_prob_E}
    (\matr{I}_{LM} \otimes \matr M - \Delta t(\matr{I}_L\otimes\matr{Q}\otimes\matr{I}_N){\vect F} - \matr{E}\otimes\matr{H})(\vect u) = \vect b
\end{align}
where the matrix $\matr{E}\in\mathbb{R}^{L\times L}$ has ones on the lower off-diagonal and zeros elsewhere, accounting for the transfer of the solution from one step to the next.

There are two fundamentally different ways to solve this system iteratively with SDC. 
We can choose either (note the $\QDmat$ instead of the $\Qmat$)
  \begin{align*} %\label{coarseP}
    \matr{P}^{\mathrm{par}}_{\vect F}({\vect u}) \coloneqq (\matr{I}_{LM}\otimes\matr{M} - \Delta t(\matr{I}_L\otimes\QDmat\otimes\matr{I}_{N}){\vect F})(\vect u)
  \end{align*}
%i.e., a block-diagonal preconditioner with SDC in the blocks (note the $\QDmat$ instead of the $\Qmat$), or
or
\begin{align*} %\label{coarseP}
    \matr{P}^{\mathrm{seq}}_{\vect F}({\vect u}) \coloneqq (\matr{I}_{LM}\otimes\matr{M} - \Delta t(\matr{I}_L\otimes\QDmat\otimes\matr{I}_{N}){\vect F} - \matr{E}\otimes\matr{{H}})(\vect u),
\end{align*}
where there are still the $\matr{H}$ matrices in the lower off-diagonal.
$\matr{P}^{\mathrm{par}}_{\vect F}$ is a parallel preconditioner, which performs SDC iterations on each step simultaneously, while $\matr{P}^{\mathrm{seq}}_{\vect F}$ propagates a single SDC iteration sequentially forward in time.

The idea of PFASST now is to couple both preconditioners in a two-level space-time full approximations scheme: the parallel $\matr{P}^{\mathrm{par}}_{\vect F}$ is used on the original problem in space and time (the ``fine'' level), while the sequential $\matr{P}^{\mathrm{seq}}_{\vect F}$ with better convergence properties is used on a coarser, cheaper level with reduced accuracy in space and/or time to reduce the impact of its sequential nature.
To create the coarse level, we reduce the number of degrees of freedom in space and choose a finite element subspace $\matr{\tilde V}^h \subset \matr{V}^h$. 
Three different transfer operations are then needed for PFASST: 
\begin{enumerate}
    \item Restriction of a coefficient-vector $\vect u_{l,m}$, representing an object in $\matr V^h$, to the representation of an object in $\matr{\tilde V}^h$,
    \item Restriction of the residual $\matr{C}_{\vect F}( \vect u) - \vect b$,
    \item Prolongation of a coefficient vector $\vect{\tilde u}_{l,m}$, representing an object in $\matr{\tilde V}_h$, to the representation of an object in $\matr{V}_h$.
\end{enumerate} 
Using Lagrange polynomials, operations 2 and 3 can be done using the canonical injection $\matr T^N\in \mathbb{R}^{N \times \tilde N}$ for prolongation and its transpose $(\matr T^N)^T \in \mathbb{R}^{\tilde N \times N}$ for restriction of the residual. 
For Operation 1, we use the matrix $\matr R^N\in\mathbb{R}^{\tilde N \times N}$ that represents the Lagrange interpolation of functions from $\matr V^h$ in $\matr{\tilde{V}}^h$.
%For Operation 1, we define the matrix $\matr R^N\in\mathbb{R}^{N \times \tilde N}$ that transfers the values of the function represented by $\vect{u}_{l,m}$ to the nodes of $\matr{\tilde V}_h$. 
By $\matr T := \matr I_{LM} \otimes \matr T^N$ and $\matr R := \matr I_{LM} \otimes \matr R^N$ we define global transfer operators.
Using the tilde symbols to indicate entities on the coarse level, one iteration of PFASST reads:%consists of the following steps: 
\begin{subequations}\label{alg_pfasst}  
  \begin{enumerate}
       \item Restrict current iterate to the coarse level: $\vect{\tilde u}^{k} =  \matr R \vect{ u}^{k}$.
       \item Compute FAS correction: $\vect\tau = \matr{\tilde C}_{\vect F} ({\vect{\tilde u}}^k) - \matr T^T \matr{ C}_{\vect F} ( {\vect u}^k)$ %\in \mathbb{R}^{L \tilde M}$.
   \item Compute $\vect{\tilde{u}}^{k+1}$ by solving $\matr{\tilde P}_{\vect F} (\vect{\tilde{u}}^{k+1} ) = (\matr{\tilde P}_{\vect F} - \matr{\tilde C}_{\vect F})({\vect{ \tilde{u}}}^{k}) +  \vect{\tilde b}+ \vect\tau$.
   \item Apply coarse grid correction: $\vect{u}^{k+\frac{1}{2}} = \vect{u}^{k} + \matr T ( \vect{\tilde{u}}^{k+1} -\matr R \vect{u}^k )$.
   \item Compute $\vect{{u}}^{k+1}$ by solving $\matr{ P}_{\vect F} ( \vect{u}^{k+1} ) = (\matr{ P}_{\vect F} - \matr{C}_{\vect F})( \vect{u}^{k+\frac{1}{2}} )  + \vect b$.
  \end{enumerate}
\end{subequations}

In contrast to the description in~\cite{bolten2017multigrid,bolten2018asymptotic}, the mass matrices are now included in $\matr{P}_{\vect F}({\vect u})$, $\matr{\tilde P}_{\vect F}(\tilde{\vect u})$ as well as in $\matr{ C}_{\vect F}$ and $\tilde{\matr{ C}}_{\vect F}$. 
This approach is preferable to others, including the naive one where the collocation problem~\eqref{eq:coll_prob} is multiplied by $\matr{M}^{-1}$ from the left. 
The collocation problem~\eqref{eq:coll_prob} then reads
\begin{align}
    \matr{\bar C}^{\operatorname{coll}}_{\vect{\bar{f}}}(\vect{u}_l) \coloneqq (\matr{I}_{M} \otimes \matr{I}_N  - \dt(\Qmat \otimes \matr I_N)\vect{\bar f} )(\vect{u}_l) = \vect{u}_{l,0}, \label{invert_coll}
\end{align}
for 
$
    \vect{\bar f}(\vect{u}_l) \coloneqq  (\bar f(u_{l,1}), ..., \bar f(u_{l,M}))^T
$
and $\bar f(u_{l,m}) = \matr M^{-1}f(u_{l,m})$. 
Similarly, the composite collocation problem then is 
\begin{align}
    \matr{\bar C}_{\vect{\bar{F}}}(\vect{u}) \coloneqq(\matr{I}_{LMN} - \Delta t(\matr{I}_L\otimes\matr{Q}\otimes\matr{I}_N){\vect{\bar{F}}} - \matr{E}\otimes\matr{N}\otimes\matr{I}_N)(\vect u) = \vect{\bar{b}}, \label{invert_comp_coll}
\end{align}
where $\vect{\bar{b}}\coloneqq (\vect{u}_{0,0}, \vect 0, ..., \vect 0)^T$ and $\vect{\bar{F}} \coloneqq \left( \vect{\bar f} ({\vect u}_{1}), \dots , \vect{\bar f} ({\vect u}_{L})  \right)^T$.
SDC and PFASST can then be derived precisely as in the literature, using the modified right-hand side $\bar f$.
Note that the inversion of the mass matrix is only necessary, if the actual residual of \eqref{invert_coll} or \eqref{invert_comp_coll} needs to be computed, which, e.g., is necessary for the FAS correction.
There, $\vect{\bar{F}}$ has to be evaluated on the fine and the coarse level, both containing the inverse of the respective mass matrix.
While seemingly attractive in terms of writing a generic code, inversion of the mass matrix can be costly and, as we will see later, convergence of PFASST is way worse in this case.
Note that the components of the residual of \eqref{invert_coll} and \eqref{invert_comp_coll} in contrast to \eqref{eq:coll_prob} and \eqref{eq:comp_coll_prob_E} are not elements of the dual space of $\matr V^h$ and therefore cannot be restricted exactly. 
In this case, the obvious choice is to use $\matr R^N$ to transfer both residual and coefficient vectors to the coarse level.

\section{Numerical results}

We now investigate numerically the convergence behavior of PFASST with finite elements in space.
In \cite{kremling2020convergence} it was shown that for a discretization with finite differences, the single-step version of PFASST (i.e.~multilevel SDC) can gain two orders of accuracy per iteration, provided very high-order transfer operators in space are used~\cite{bolten2017multigrid,kremling2020convergence}.
We will now show numerically that with finite elements in space, this is no longer necessary.

We use the following nonlinear differential equation
\begin{align}
u_t = \Delta u + u^2 (1 - u) \quad \text{ on } [0,2]\times [-20,20]. \label{eq:flame}
\end{align}
In all simulations, we use $4$ Gauss--Raudau nodes to discretize a single time-step. 
In the following we use SDC for serial time-step calculations and PFASST to calculate $4$ time-steps simultaneously.
The spatial domain $[-20,20]$ is discretized using Lagrange finite elements of either order $1$ or order $3$.
We use the initial value $u(x,0) = (1+(\sqrt{2}-1)e^{-\sqrt{6}/6x})^{-2}$ as the initial guess for the iteration.

For the first test case we use a third-order Lagrange basis to approximate the solution.
We coarsen the problem in space by restricting to a second-order Lagrange basis.
Figure \ref{hiSDC} shows the results for SDC and PFASST.
They show the absolute error of the method in the infinity norm for different time-step sizes, in relation to a reference solution calculated with a much smaller $\dt$ and SDC. 
While SDC gains one order per iteration as expected, PFASST can gain up to two orders per iteration, at least after some initial iterations have been performed.
This ``burn-in'' phase causes a loss of parallel efficiency when actual speedup is measured.
After this phase, however, PFASST shows ideal convergence behavior gaining two orders of accuracy per iteration.
There is not yet a theoretical explanation for neither the ``burn-in'' nor the ``ideal'' phase.

For a second test case we use a first-order Lagrange basis and coarsen the problem in space by doubling the element size.
Figure \ref{linearSDC} shows the results for SDC and PFASST.
In the same way as in the high-order example before, SDC gains one order or accuracy per iteration, while PFASST can gain up to two orders after a few initial iterations.
Note that in the case of a finite difference discretization, the order of the interpolation is crucial to obtain two orders per iteration~\cite{bolten2017multigrid,kremling2020convergence}.
The usage of $\matr{T}^N$ as exact interpolation for nested finite element spaces removes this, so far, persistent and irritating limitation.

\begin{figure}[t]
\begin{tikzpicture}[scale=0.64]	\begin{loglogaxis}[legend pos=south west,xlabel= $\Delta t$,	ylabel= error,	x dir=reverse, label style={font=\large}, , xticklabels={$0.5$,$0.25$,$0.125$,$0.0625$, $0.03125$},xtick={0.5,0.25,0.125,0.0625, 0.03125}, ymax=1,ymin=10e-14	]  
	\addplot [red, mark=*] plot coordinates {
		(0.5, 0.0118369)
		(0.25, 0.00583352)
		(0.125, 0.00289607)
		(0.0625, 0.00144291)
		(0.03125, 0.000720183)
	};
	\addplot [domain= 0.03125:0.5, dotted, black, thick] {(x* 0.0118369/0.5)} node[below] {order 1};
	\addplot [brown, mark=square*] plot coordinates {
		(0.5, 0.00080651)
		(0.25, 0.00018457)
		(0.125, 4.42428e-05)
		(0.0625, 1.08348e-05)
		(0.03125, 2.68178e-06)
	};
	\addplot [domain= 0.03125:0.5, dotted, black, thick] {(x^2*2* 0.00080651/0.5)} node[below] {order 2};
	\addplot [green, mark=otimes*] plot coordinates {
		(0.5, 6.18404e-05)
		(0.25, 6.37561e-06)
		(0.125, 7.28366e-07)
		(0.0625, 8.71624e-08)
		(0.03125, 1.06643e-08)
	};
	\addplot [domain= 0.03125:0.5, dotted, black, thick] {(x^3*3* 6.18404e-05/0.5)} node[below] {order 3};
	\addplot [gray, mark=diamond*] plot coordinates {
		(0.5, 5.15203e-06)
		(0.25, 2.35574e-07)
		(0.125, 1.26734e-08)
		(0.0625, 7.35936e-10)
		(0.03125, 4.42273e-11)
	};
	\addplot [domain= 0.03125:0.5, dotted, black, thick] {(x^4*4* 5.15203e-06/0.5)} node[below] {order 4};
	\addplot [blue, mark=star] plot coordinates {
		(0.5, 4.67086e-07)
		(0.25, 9.72257e-09)
		(0.125, 2.52354e-10)
		(0.0625, 7.3635e-12)
		(0.03125, 3.54959e-13)
	};
	\addplot [domain= 0.03125:0.5, dotted, black, thick] {(x^5*5* 4.67086e-07/0.5)} node[below] {order 5};
	\legend{k=1,,k=2,,k=3,,k=4,,k=5}	\end{loglogaxis} \end{tikzpicture}
	\begin{tikzpicture}[scale=0.64]	\begin{loglogaxis}[legend pos=outer north east,xlabel= $\Delta t$,	x dir=reverse, label style={font=\large}, , xticklabels={$0.5$,$0.25$,$0.125$,$0.0625$, $0.03125$},xtick={0.5,0.25,0.125,0.0625, 0.03125}, , ymax=1,ymin=10e-14	]  
	\addplot [red, mark=*] plot coordinates {
		(0.5, 0.112573)
		(0.25, 0.117488)
		(0.125, 0.12012)
		(0.0625, 0.121477)
		(0.03125, 0.122155)
	};
	\addplot [brown, mark=square*] plot coordinates {
		(0.5, 0.122925)
		(0.25, 0.122905)
		(0.125, 0.122871)
		(0.0625, 0.12286)
		(0.03125, 0.122857)
	};
	\addplot [green, mark=otimes*] plot coordinates {
		(0.5, 0.00510645)
		(0.25, 0.00373913)
		(0.125, 0.00213376)
		(0.0625, 0.00112915)
		(0.03125, 0.000579769)
	};
	\addplot [domain = 0.03125:0.5, dotted, black, thick] {(x^1*1*0.00510645/0.5)} node[below] {order 1};
	\addplot [gray, mark=diamond*] plot coordinates {
		(0.5, 0.000222438)
		(0.25, 2.87173e-05)
		(0.125, 3.52793e-06)
		(0.0625, 4.34549e-07)
		(0.03125, 5.38507e-08)
	};
	\addplot [domain = 0.03125:0.5, dotted, black, thick] {(x^3*3*0.000222438/0.5)} node[below] {order 3};
	\addplot [blue, mark=star] plot coordinates {
		(0.5, 6.32669e-06)
		(0.25, 1.5478e-07)
		(0.125, 9.72825e-09)
		(0.0625, 4.86207e-09)
		(0.03125, 8.38465e-10)
	};
	\addplot [domain = 0.03125:0.5, dotted, black, thick] {(x^5*5*6.32669e-06/0.5)} node[below] {order 5};
	\end{loglogaxis} \end{tikzpicture} \caption{SDC (left) and PFASST (right) errors for different $\dt$ and number of iterations $k$, $128$ spatial elements, order $3$} \label{hiSDC}\end{figure}
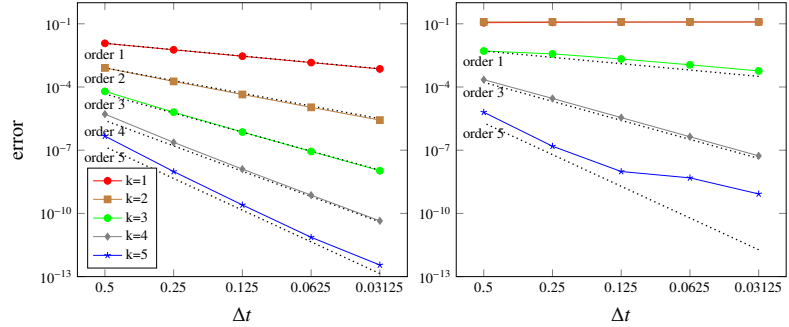
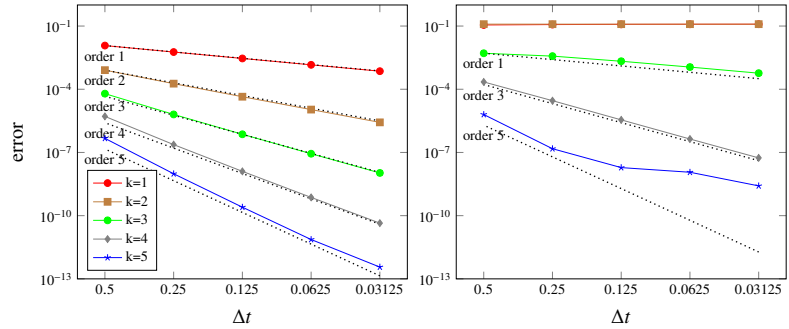
\begin{figure}[t] 
	\begin{tikzpicture}[scale=0.64]	\begin{loglogaxis}[legend pos=south west,xlabel= $\Delta t$,	ylabel= error,	x dir=reverse, label style={font=\large}, , xticklabels={$0.5$,$0.25$,$0.125$,$0.0625$, $0.03125$},xtick={0.5,0.25,0.125,0.0625, 0.03125}, ymax=1,ymin=10e-14	]  
	\addplot [red, mark=*] plot coordinates {
		(0.5, 0.0118389)
		(0.25, 0.00583435)
		(0.125, 0.00289646)
		(0.0625, 0.0014431)
		(0.03125, 0.000720277)
	};
	\addplot [domain= 0.03125:0.5, dotted, black, thick] {(x*0.0118389/0.5)} node[below] {order 1};
	\addplot [brown, mark=square*] plot coordinates {
		(0.5, 0.00080651)
		(0.25, 0.000184572)
		(0.125, 4.42463e-05)
		(0.0625, 1.08381e-05)
		(0.03125, 2.68224e-06)
	};
	\addplot [domain= 0.03125:0.5, dotted, black, thick] {(x^2*2*0.00080651/0.5)} node[below] {order 2};
	\addplot [green, mark=otimes*] plot coordinates {
		(0.5, 6.18517e-05)
		(0.25, 6.37619e-06)
		(0.125, 7.28373e-07)
		(0.0625, 8.71693e-08)
		(0.03125, 1.06663e-08)
	};
	\addplot [domain= 0.03125:0.5, dotted, black, thick] {(x^3*3*6.18517e-05/0.5)} node[below] {order 3};
	\addplot [gray, mark=diamond*] plot coordinates {
		(0.5, 5.15238e-06)
		(0.25, 2.35568e-07)
		(0.125, 1.26728e-08)
		(0.0625, 7.35906e-10)
		(0.03125, 4.42195e-11)
	};
	\addplot [domain= 0.03125:0.5, dotted, black, thick] {(x^4*4*5.15238e-06/0.5)} node[below] {order 4};
	\addplot [blue, mark=star] plot coordinates {
		(0.5, 4.67103e-07)
		(0.25, 9.72289e-09)
		(0.125, 2.5239e-10)
		(0.0625, 7.37163e-12)
		(0.03125, 3.6237e-13)
	};
	\addplot [domain= 0.03125:0.5, dotted, black, thick] {(x^5*5*4.67103e-07/0.5)} node[below] {order 5};
	\legend{k=1,,k=2,,k=3,,k=4,,k=5,}	\end{loglogaxis} \end{tikzpicture} 
	\begin{tikzpicture}[scale=0.64]	\begin{loglogaxis}[legend pos=outer north east,xlabel= $\Delta t$,	x dir=reverse, label style={font=\large}, , xticklabels={$0.5$,$0.25$,$0.125$,$0.0625$, $0.03125$},xtick={0.5,0.25,0.125,0.0625, 0.03125}, , ymax=1,ymin=10e-14	]  
	\addplot [red, mark=*] plot coordinates {
		(0.5, 0.112587)
		(0.25, 0.1175)
		(0.125, 0.120136)
		(0.0625, 0.121493)
		(0.03125, 0.12218)
	};
	\addplot [brown, mark=square*] plot coordinates {
		(0.5, 0.122927)
		(0.25, 0.122911)
		(0.125, 0.122878)
		(0.0625, 0.122868)
		(0.03125, 0.122865)
	};
	\addplot [green, mark=otimes*] plot coordinates {
		(0.5, 0.00510769)
		(0.25, 0.00373968)
		(0.125, 0.00213391)
		(0.0625, 0.00112933)
		(0.03125, 0.000579787)
	};
	\addplot [domain = 0.03125:0.5, dotted, black, thick] {x*0.00510769/0.5} node[below] {order 1};
	\addplot [gray, mark=diamond*] plot coordinates {
		(0.5, 0.000222391)
		(0.25, 2.86896e-05)
		(0.125, 3.52339e-06)
		(0.0625, 4.34882e-07)
		(0.03125, 5.52076e-08)
	};
	\addplot [domain = 0.03125:0.5, dotted, black, thick] {(x^3*3*0.000222391/0.5)} node[below] {order 3};
	\addplot [blue, mark=star] plot coordinates {
		(0.5, 6.31837e-06)
		(0.25, 1.50525e-07)
		(0.125, 1.91034e-08)
		(0.0625, 1.14888e-08)
		(0.03125, 2.59101e-09)
	};
	\addplot [domain = 0.03125:0.5, dotted, black, thick] {(x^5*5*6.31837e-06/0.5)} node[below] {order 5};
	\end{loglogaxis} \end{tikzpicture} 
	\caption{SDC (left) and PFASST (right) errors for different $\dt$ and number of iterations $k$, $512$ spatial elements, order $1$}
	\label{linearSDC}
	\end{figure}
	\begin{figure}[h!] \begin{tikzpicture}[scale=0.64]	\begin{loglogaxis}[legend pos=south west,xlabel= $\Delta t$,	ylabel= error,	x dir=reverse, label style={font=\large}, , xticklabels={$0.5$,$0.25$,$0.125$,$0.0625$, $0.03125$},xtick={0.5,0.25,0.125,0.0625, 0.03125}, , ymax=1,ymin=10e-14	]  
\addplot [red, mark=*] plot coordinates {
	(0.5, 0.0118389)
	(0.25, 0.00583435)
	(0.125, 0.00289646)
	(0.0625, 0.0014431)
	(0.03125, 0.000720277)
};
	\addplot [domain = 0.03125:0.5, dotted, black, thick] {(x^1*1*0.0118389/0.5)} node[below] {order 1};
\addplot  [brown, mark=square*] plot coordinates {
	(0.5, 0.00080651)
	(0.25, 0.000184572)
	(0.125, 4.42463e-05)
	(0.0625, 1.08381e-05)
	(0.03125, 2.68224e-06)
};
	\addplot [domain = 0.03125:0.5, dotted, black, thick] {(x^2*2*0.00080651/0.5)} node[below] {order 2};
\addplot  [green, mark=otimes*]  plot coordinates {
	(0.5, 6.18517e-05)
	(0.25, 6.37619e-06)
	(0.125, 7.28373e-07)
	(0.0625, 8.71692e-08)
	(0.03125, 1.06662e-08)
};
	\addplot [domain = 0.03125:0.5, dotted, black, thick] {(x^3*3*6.18517e-05/0.5)} node[below] {order 3};
\addplot  [gray, mark=diamond*] plot coordinates {
	(0.5, 5.15238e-06)
	(0.25, 2.35569e-07)
	(0.125, 1.2673e-08)
	(0.0625, 7.36055e-10)
	(0.03125, 4.43684e-11)
};
	\addplot [domain = 0.03125:0.5, dotted, black, thick] {(x^4*4*5.15238e-06/0.5)} node[below] {order 4};
\addplot [blue, mark=star] plot coordinates {
	(0.5, 4.67103e-07)
	(0.25, 9.72274e-09)
	(0.125, 2.52243e-10)
	(0.0625, 7.22566e-12)
	(0.03125, 2.16459e-13)
};
	\addplot [domain = 0.03125:0.5, dotted, black, thick] {(x^5*5*4.67103e-07/0.5)} node[align=right] {order 5};
	\legend{k=1,,k=2,,k=3,,k=4,,k=5} \end{loglogaxis}  \end{tikzpicture} 
%\caption{SDC (left) and PFASST (right) convergence for different number of iterations, $4$  nodes per time-step, $512$ spatial elements, order $1$}
%\label{linearSDC}
	\begin{tikzpicture}[scale=0.64]	\begin{loglogaxis}[legend pos=outer north east,xlabel= $\Delta t$,	x dir=reverse, label style={font=\large},  xticklabels={$0.5$,$0.25$,$0.125$,$0.0625$, $0.03125$},xtick={0.5,0.25,0.125,0.0625, 0.03125}, , ymax=1,ymin=10e-14	]  
\addplot  [red, mark=*] plot coordinates {
	(0.5, 0.118809)
	(0.25, 0.121014)
	(0.125, 0.12199)
	(0.0625, 0.122409)
	(0.03125, 0.122633)
}; %[red, mark=*]  [brown, mark=square*] [green, mark=otimes*]  [gray, mark=diamond*]  [blue, mark=star]
\addplot [brown, mark=square*] plot coordinates {
	(0.5, 0.123132)
	(0.25, 0.123176)
	(0.125, 0.12296)
	(0.0625, 0.122886)
	(0.03125, 0.122867)
};
\addplot [green, mark=otimes*] plot coordinates {
	(0.5, 0.0025012)
	(0.25, 0.00067517)
	(0.125, 0.000323288)
	(0.0625, 0.000210587)
	(0.03125, 0.000121386)
};
\addplot [domain= 0.03125:0.5, dotted, black, thick] {(x*0.0025012/0.5)} node[below] {order 1};
\addplot [gray, mark=diamond*] plot coordinates {
	(0.5, 0.000228126)
	(0.25, 0.000125566)
	(0.125, 7.6066e-05)
	(0.0625, 4.07373e-05)
	(0.03125, 2.00234e-05)
};
%	\addplot [domain= 0.03125:0.5, dotted, black, thick] {(x*0.000228126/0.5)} node[below] {};
\addplot [blue, mark=star] plot coordinates {
	(0.5, 0.000143917)
	(0.25, 5.47294e-05)
	(0.125, 1.47057e-05)
	(0.0625, 6.49487e-06)
	(0.03125, 3.51247e-06)
};
%	\addplot [domain= 0.03125:0.5, dotted, black, thick] {(x* 0.000143917/0.5)} node[below] {};
	\end{loglogaxis} \end{tikzpicture}
	\caption{Naive approach with inverted mass matrix: SDC (left) and PFASST (right) errors for different $\dt$ and number of iterations $k$, $512$ spatial elements, order~$1$}
	\label{MinvPFASST}
\end{figure}
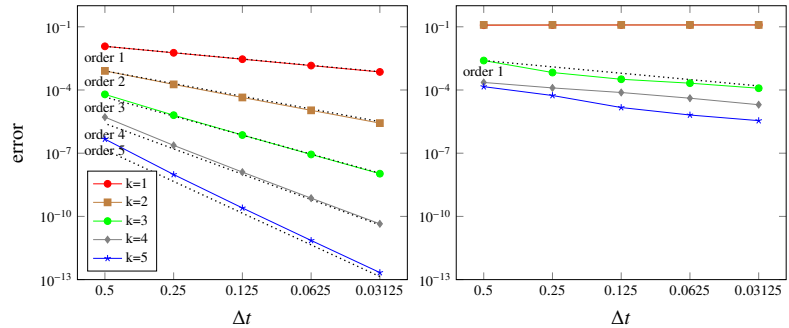

Finally, Figure~\ref{MinvPFASST} shows SDC and PFASST applied to the (composite) collocation problem~\eqref{invert_comp_coll} with inverted mass matrix and a first-order Lagrange basis.
SDC behaves exactly as before, while PFASST fails to show any reasonable convergence. 
In particular, increasing the number of iterations does not increase the order of accuracy beyond $1$.

The advantage of using finite elements together with PFASST in the way we demonstrated here is not yet analyzed analytically. 
We intend to address this in a follow-up work.
Also, the important fact that two orders of accuracy per iteration is possible even with a low-order spatial interpolation does not have a theoretical explanation. 
A corresponding analysis is work in progress. 

% \input{sections/conclusions}

%\bibliographystyle{abbrv}
%\bibliography{refs}
%\bibliography{}
\end{document}